\documentclass[12pt]{amsart}
\usepackage{amsmath,amssymb,amsthm}
\usepackage[]{latexsym}
\newtheorem{defn0}{Definition}[section]
\newtheorem{prop0}[defn0]{Proposition}
\newtheorem{thm0}[defn0]{Theorem}
\newtheorem{lemma0}[defn0]{Lemma}
\newtheorem{corollary0}[defn0]{Corollary}
\newtheorem{example0}[defn0]{Example}
\newtheorem{remark0}[defn0]{Remark}
\newtheorem{conjecture0}[defn0]{Conjecture}
\newtheorem{notation0}[defn0]{Notation}

\newenvironment{definition}{\begin{defn0}\rm}{\end{defn0}}
\newenvironment{proposition}{\begin{prop0}}{\end{prop0}}
\newenvironment{theorem}{\begin{thm0}}{\end{thm0}}
\newenvironment{lemma}{\begin{lemma0}}{\end{lemma0}}
\newenvironment{corollary}{\begin{corollary0}}{\end{corollary0}}
\newenvironment{example}{\begin{example0}\rm}{\end{example0}}
\newenvironment{remark}{\begin{remark0}\rm}{\end{remark0}}

\newcommand{\Gal}{{\mathrm {Gal }}}
\newcommand{\n}{\mathrm{n}}
\newcommand{\tr}{\mathrm{tr}}

\newcommand{\Pic}{\mathrm{Pic}}
\newcommand{\Norm}{\mathrm{N}}
\newcommand{\Aut}{\mathrm{Aut}}

\newcommand{\End}{{\mathrm{End}}}

\newcommand{\disc}{{\mathrm {disc }}}

\newcommand{\M}{{\mathrm M}}

\newcommand{\Z}{{\mathbb Z}}
\newcommand{\Q}{{\mathbb Q}}
\newcommand{\C}{{\mathbb C}}
\newcommand{\R}{{\mathbb R}}

\newcommand{\HH}{{\mathbb {H}}}

\newcommand{\cL}{{\mathcal L}}

\newcommand{\cH}{{\mathcal H}}

\newcommand{\cN}{{\mathcal N}}

\newcommand{\cO}{{\mathcal O}}
\newcommand{\cT}{{\mathcal T}}

\newcommand{\ra}{\rightarrow}
\newcommand{\lra}{\longrightarrow}

\newcommand{\qbar}{{\overline \Q}}

\newcommand{\om}{{\omega}}
\newcommand{\sg}{{\sigma}}

\include{thebibliography}

\begin{document}

\title{ The field of moduli of quaternionic\\
multiplication on abelian varieties }

\author{Victor Rotger}
\footnote{Partially supported by a grant FPI from the Ministerio
de Ciencia y Tecnolog\'{\i}a and by BFM2000-0627. }

\address{Escola Universit\`aria Polit\`ecnica de
Vilanova i la Geltr\'u, Av. V\'{\i}ctor Balaguer s/n, E-08800
Vilanova i la Geltr\'u, Spain}

\email{vrotger@mat.upc.es}



\begin{abstract}

We consider principally polarized abelian varieties with
quaternionic multiplication over number fields and we study the
field of moduli of their endomorphisms in relation to the set of
rational points on suitable Shimura varieties. \\ \\ {\smaller
Published in {\em Intern.\ J.\ Math.\ M.\ Sc.\ } {\bf 52} (2004),
2795-2808.}
\end{abstract}

\maketitle

\section{Introduction}

Let $\bar \Q $ be a fixed algebraic closure of the field $\Q $ of
rational numbers and let $(A, \cL )/\bar \Q $ be a polarized abelian variety. The
field of moduli of $(A, \cL )/\bar \Q $ is the minimal number field
$k_{A, \cL }\subset \qbar $ such that $(A, \cL )$ is isomorphic (over $\qbar $) to
its Galois conjugate $(A^{\sigma }, \cL ^{\sg })$, for all $\sigma \in
\Gal (\qbar /k_{A, \cL })$.

The field of moduli $k_{A, \cL }$ is an essential arithmetic
invariant of the $\qbar $-isomorphism class of $(A, \cL )$. It is
contained in all possible fields of definition of $(A, \cL )$ and,
unless $(A, \cL )$ admits a rational model over $k_{A, \cL }$
itself, there is not a unique minimal field of definition for $(A,
\cL )$. In this regard, we have the following theorem of Shimura.

\begin{theorem}[\cite{Sh2}]\label{Shh}
A generic principally polarized abelian variety of odd dimension
admits a model over its field of moduli. For a generic principally
polarized abelian variety of even dimension, the field of moduli
is not a field of definition.

\end{theorem}

Let $\End (A)=\End _{\qbar }(A)$ denote the ring of endomorphisms
of $A$. It is well known that $\End (A)=\Z $ for a generic polarized
abelian variety $(A, \cL  )$. However, due to Albert's classification of involuting
division algebras (\cite{Mu}) and the work of Shimura (\cite{Sh1}, there are
other rings that can occur as the
endomorphism ring of an abelian variety. Namely, if $A/\qbar $ is
simple, $\End (A)$ is an order in either a totally
real number field $F$ of degree $[F:\Q ]\mid \dim (A)$, a totally
indefinite quaternion algebra $B$ over a totally real number field
$F$ of degree $2 [F:\Q ]\mid \dim (A)$, a totally definite
quaternion algebra $B$ over a totally real number field $F$ of
degree $2 [F:\Q ]\mid \dim (A)$ or a division algebra over a
CM-field.

Let us recall that a quaternion algebra $B$ over a totally real
field $F$ is called totally indefinite if $B\otimes _{\Q }\R
\simeq \M _2(\R ) \oplus ... \oplus \M _2(\R )$ and totally
definite if $B\otimes _{\Q }\R \simeq \HH \oplus ... \oplus \HH $,
where $\HH = (\frac {-1, -1}{\R })$ denotes the skew-field of real
Hamilton quaternions.

\begin{definition}

Let $(A, \cL )/\qbar $ be a polarized abelian variety and let
$S\subseteq \End (A)$ be a subring of $\End (A)$. The field of moduli of
$S$ is the minimal number
field $k_S\supseteq k_{A, \cL }$ such that, for any $\sigma \in \Gal
(\qbar /k_S)$,
there is an isomorphism $\varphi _{\sigma } /\qbar :A\ra A^{\sigma
}$, $\varphi _{\sg }^*(\cL ^{\sg }) = \cL $, of polarized abelian varieties
that induces commutative diagrams

$$
\quad A\lra A^{\sg }
$$
$$
\quad \beta \downarrow \quad \quad \quad \downarrow \beta ^{\sg }
$$
$$
\quad A\lra A^{\sg }
$$
for any $\beta \in S$.
\end{definition}

We remark that, as a consequence of the very basic definitions,
the field of moduli of the multiplication-by-$n$ endomorphisms on
$A$ is exactly $k_{\Z } = k_{A, \cL }$. But in the case that $\End
(A)\varsupsetneq \Z $, little is known on the chain of Galois
extensions $k_{\End (A)}\supseteq k_S \supseteq k_{A, \cL }$.

The main aim of this article is to study the field of moduli
of totally indefinite quaternionic multiplication on an abelian variety.
In relation to Shimura's Theorem \ref{Shh}, we remark that the
dimension of an abelian variety whose endomorphism ring contains a quaternion
order is always even.

We state our main result in the next section. As we will show in
Section \ref{ff}, it is a consequence of the results obtained in
\cite{Ro2}, \cite{Ro3} on certain modular forgetful morphisms
between certain Shimura varieties, Hilbert modular varieties and
the moduli spaces of principally polarized abelian varieties.

In Section \ref{defi}, we particularize our results to abelian
surfaces. We use our results together with those of Mestre \cite{Me} and Jordan
\cite{Jo} to compare the field of moduli and field of definition of the
quaternionic multiplication on an abelian surface.

In an appendix to this paper, we discuss a question on the
arithmetic of quaternion algebras that naturally arises from our
considerations and which is also related to recent work by
Chinburg and Friedman \cite{ChFr1}, \cite{ChFr2}.

A cryptographical application of the results in the appendix has
been derived in \cite{GR} by Galbraith and the author.

\section{Main result}\label{mm}

Let $F$ be a totally real number field $F$ of degree $[F:\Q ]=n$
and let $R_F$ denote its ring of integers. We shall let $F^*_+$
denote the subgroup of totally positive elements of $F^*$. For any
finite field extension $L/F$, let $R_L$ denote the ring of
integers of $L$ and let $\Omega _{odd}(L)= \{ \xi \in R_L, \xi ^{f}=1,
f\mbox{ odd }\} $ denote the set of primitive roots of unity of
odd order in $L$. We let $\omega _{odd}(L) = |\Omega (L)|$.

Let $B$ be a totally indefinite quaternion algebra over $F$ and
let $\cO$ be a maximal order in $B$.

\begin{definition}

An abelian variety $A/k$ over an algebraically closed field $k$
has quaternionic multiplication by $\cO$ if $\End (A)\simeq \cO
_B$ and $\dim (A) = 2 n$.

\end{definition}

\begin{proposition}\cite{Ro1}
Let $(A, \cL )$ be a principally polarized abelian variety with
quaternionic multiplication by $\cO$ over $\qbar $. Then the
discriminant ideal $\disc (B)$ of $B$ is principal and generated
by a totally positive element $D\in F^*_+$.

\end{proposition}

As in \cite{Ro2}, \cite{Ro3} we say that a quaternion algebra $B$
over $F$ of totally positive principal discriminant $\disc (B)\in
F^*_+$ is {\em twisting} if $B\simeq (\frac {-D, m}{F})$ for some
$m\in F^*_+$ supported at the prime ideals $\wp \mid D$ of $F$.
Let $C_2$ denote the cyclic group of two elements. The main result
of this article is the following.

\begin{theorem}
\label{main} Let $(A, \cL )$ be a principally polarized abelian
variety with quaternionic multiplication by $\cO$ over $\qbar $
and let $\disc (B) =D$ for some $D\in F^*_+$. Let $\omega = \omega
(F(\sqrt {-D})$.
\begin{enumerate}
\item[(i)]
If $B$ is not twisting, then

\begin{itemize}
\item For any totally real quadratic order $S\subset \cO$ over
$R_F$, $k_{\cO }=k_S$.

\item $\Gal (k_{\cO }/k_{R_F})\subseteq C_2^{\omega _{odd}}$.
\end{itemize}

\item[(ii)]
If $B$ is twisting, then

\begin{itemize}
\item For any totally real quadratic order $S\subset \cO$,
$\Gal (k_{\cO }/k_{S})\subseteq C_2$.

\item $\Gal (k_{\cO }/k_{R_F})\subseteq C_2^{2 \omega _{odd}}$.

\end{itemize}

\end{enumerate}

\end{theorem}

As we state more precisely in Section \ref{ff}, Theorem \ref{main}
admits several refinements.

\section{Proof of Theorem \ref{main}: Shimura varieties and forgetful
maps}\label{ff}

Let $B$ be a totally indefinite quaternion division algebra over a
totally real number field $F$ and assume that $\disc (B) = (D)$ for
some $D\in F^*_+$. Let $\cO$ be a maximal order in $B$ and fix an
arbitrary quaternion $\mu \in \cO $ satisfying
$\mu ^2+D=0$. Its existence is
guaranteed by Eichler's theory on optimal embeddings (\cite{Vi}) and it
generates a CM-field $F(\mu )\simeq F(\sqrt
{-D})$ over $F$ embedded in $B$. We will refer to the pair $(\cO, \mu )$ as
a {\em principally polarized
order}.

Attached to $(\cO,\mu )$, we can consider a Shimura
variety $X_{\cO, \mu }/\Q $ that solves the coarse moduli
problem of classifying triplets $(A, \iota , \cL )$ over $\Q $
where:

\begin{enumerate}

\item[(i)]
$(A, \cL )$ is a principally polarized abelian variety and
\item[(ii)]
$\iota : \cO\hookrightarrow \End (A)$ is a monomorphism of rings
satisfying that $\iota (\beta )^* = \iota (\mu ^{-1} \bar {\beta }\mu )$ for
all $\beta \in \cO$ and where $^*$ denotes the Rosati involution with
respect to $\cL $.

\end{enumerate}

Attached to the maximal order $\cO$ there is also the
{\em Atkin-Lehner group}

$$
W=\Norm _{B^*}(\cO)/F^* \cO^*.
$$

Let $C_2$ denote the cyclic group of two elements. The group $W$ is isomorphic
to $C_2\times \stackrel{2 r}{...} \times C_2$,
where $2 r = \sharp \{ \wp \mid \disc (B)\}$ is the
number of ramifying prime ideals of $B$ (cf.\,\cite{Vi}, \cite{Ro2}).

Let $B^*_+$ be the group of
invertible quaternions of totally positive reduced norm. The
{\em positive Atkin-Lehner group} is $W^1=\Norm _{B^*_+}(\cO )/F^*
\cO ^1$, where $\cO ^1=\{ \gamma \in \cO , \n (\gamma )=1\} $
denotes the group of units of $\cO $ of reduced norm $1$.

As it was shown in \cite{Ro2}, the group $W^1$ is a subgroup of the
automorphism group $\Aut _{\Q }(X_{\cO, \mu })$ of the Shimura
variety $X_{\cO, \mu }$.

We have

$$
W^1\simeq C_2^s \mbox{ for } 2 r\leq s\leq n + 2r -1.
$$

The first inequality holds because there is a
natural map $W^1\twoheadrightarrow W$ which is an epimorphism of groups
due to indefiniteness of $B$ and the norm theorem for maximal orders
(see \cite{Ro2}). The second inequality is a consequence of
Dirichlet's unit theorem and it is actually an equality if the narrow class number
of $F$ is $h_+(F)=1$, as is the case of $F=\Q $.

We now introduce the notion of {\em twists} of a polarized order
$(\cO,\mu )$.

\begin{definition}

Let $(\cO, \mu )$ be a principally polarized
maximal order in a totally indefinite quaternion algebra $B$
of discriminant $\disc (B)=(D)$, $D\in F^*_+$.

A twist of $(\cO,\mu )$ is an element $\chi \in \cO \cap
\Norm _{B^*}(\cO)$ such that $\chi ^2+\n (\chi )=0, \mu \chi =-\chi \mu
$ and therefore

$$
B = F+F\mu +F\chi +F \mu \chi = (\dfrac {-D, -\n (\chi )}{F}).
$$

For any subring $S\subset \cO $, we say that $\chi $ is a twist of
$(\cO , \mu )$ in $S$ if $\chi \in S$.
\end{definition}

We say that $(\cO, \mu )$ is
{\em twisting} if it admits some twist in $\cO $ and that a
quaternion algebra $B$ is twisting if it contains a twisting polarized
maximal order. This agrees with our terminology in the preceding
section.

\begin{definition}
\label{twist}

A twisting involution $\om \in W^1$ is an Atkin-Lehner involution such that
$[\om ]=[\chi ]\in
W$ is represented by a twist $\chi $ of $(\cO , \mu )$. It is a twisting
involution in $S\subseteq \cO $
if it can be represented by a twist $\chi \in S$.

Let $V_0(S)$ denote the subgroup of $W^1$ generated by the
twisting involutions of $(\cO , \mu )$ in $S$; we will simply
write $V_0$ for $V_0(\cO )$.

\end{definition}

Let us remark that, since $B$ is totally indefinite, no $\chi \in B^*_+$
can be a twist of $(\cO , \mu )$ because a necessary condition for
$B\simeq (\frac {-D, -\n (\chi )}{F})$ is that $\n (\chi )$ be
totally negative. In fact, twisting involutions $\om \in W^1$ are
always represented by twists $\chi \in B^*_-$ of totally negative reduced
norm.

Note also that a necessary and sufficient condition for $B$ to be
twisting is that $B\simeq (\frac {-D, m}{F})$ for some element
$m\in F^*_+$ supported at the prime ideals  $\wp \mid D$ (that is,
$v_{\wp }(m)\not = 0$ only if $\wp \mid D$).

For a polarized order $(\cO , \mu )$, let $R_{\mu }= F(\mu )\cap
\cO $ be the order in the CM-field $F(\mu )\simeq F(\sqrt {-D})$
that optimally embeds in $\cO $. Note that, since $\mu \in \cO $,
$R_{\mu }\supseteq R_F[\sqrt {-D}]$. We let $\Omega = \Omega
(R_{\mu })= \{ \xi \in R_{\mu }, \xi ^f=1, f\geq 1\}$ denote the
finite group of roots of unity in $R_{\mu }$ and $\Omega _{odd}=
\{ \xi \in R_{\mu }, \xi ^{f}=1, f\mbox{ odd }\} $ the subgroup of
primitive roots of unity of odd order. Their respective
cardinalities will be denoted by $\omega = \omega (R_{\mu })$ and
$\omega _{odd} = \omega _{odd}(R_{\mu })$.

\begin{definition}
\label{stable}

The {\em stable group} $W_0=U_0\cdot V_0$ associated to
$(\cO , \mu )$ is the subgroup of
$W^1$ generated by $U_0 = \Norm _{F(\mu )^*}(\cO )/F^* \cdot \Omega $ and
the group of twisting involutions $V_0$.

\end{definition}

Note that $U_0$ is indeed a subgroup of $W^1$ because $\Omega =F(\mu )
\cap \cO ^1$.

The motivation for
introducing the Shimura variety $X_{\cO , \mu }$ and the above
Atkin-Lehner groups in this note is that it gives a
modular interpretation of the field of moduli $k_{\cO }$ of the
quaternionic multiplication on $A$: $k_{\cO }= \Q (P)$
is the extension over $\Q $
generated by the coordinates of the point $P=[A, \iota , \cL ]$ on
Shimura's canonical model $X_{\cO , \mu }/\Q $ that represents the
$\qbar $-isomorphism class of the triplet.

A similar construction holds for the totally real subalgebras of $B$. Indeed,
let $L\subset B$ be a totally real quadratic extension of $F$ embedded in $B$.
Then $S = L\cap \cO $ is an order of $L$ over $R_F$ which is optimally
embedded in $\cO $. Identifying $S$ with a subring of the ring
of endomorphisms of $A$, we again have that the field of moduli $k_S$ is the
extension $\Q (P_{|S})$ of $\Q $ generated by the coordinates of
the point $P_{|S} = [A, \iota _{|S}, \cL ]$ on the
Hilbert-Blumenthal variety $\cH _S/\Q $ that solves the coarse moduli
problem of classifying abelian varieties of dimension $2 n$ with
multiplication by $S$.

Along the same lines, the field of moduli $k_{R_F}$ of the central
endomorphisms of $A$ is the extension $\Q (P_{|{R_F}})$ of $\Q $
generated by the coordinates of the point $P_{|{R_F}} = [A, \iota
_{|{R_F}}, \cL ]$ on the Hilbert-Blumenthal variety $\cH _F/\Q $
which solves the coarse moduli problem of classifying abelian
varieties of dimension $2 n$ with multiplication by $R_F$.

%Note that there are infinitely many choices of totally real
%quadratic orders $S$ in $\cO $ over $R_F$.

The tool for studying the Galois extensions $k_{\cO }/k_{S}/k_{R_F}$ is
provided by the forgetful modular maps

$$
\begin{matrix}
\pi _F: & X_{\cO , \mu } & \stackrel {\pi _S}{\longrightarrow }&
\cH _S & \longrightarrow  & \cH _F \\
&  P  &\mapsto & P_{|S} & \mapsto & P_{|R_F}.\\
\end{matrix}
$$

It was shown in \cite{Ro2} that the morphisms $\pi _F$ and $\pi _S$
have finite fibres. Furthermore, it was proved in \cite{Ro2} that:

\begin{enumerate}

\item[(i)] There is a birational equivalence $b_S: X_{\cO ,
\mu }/V_0(S)\stackrel{\sim }{\ra }\pi _{S, \varphi }(X_{\cO ,
\mu })$ and a commutative diagram
$$
\begin{matrix}
\pi _S: & X_{\cO , \mu } &           & \longrightarrow &       & \cH _S, \\
        &                & \stackrel{\pi _S}{\searrow  }&                 &  \stackrel{b_S}{\nearrow }     & \\
        &                &           & X_{\cO ,\mu }/V_0(S) &     &
\end{matrix}
$$
where $p_S: X_{\cO , \mu }\ra X_{\cO , \mu }/V_0(S)$ is the natural
projection.

\item[(ii)] There is a birational equivalence $b_S: X_{\cO ,
\mu }/V_0(S)\stackrel{\sim }{\ra }\pi _{S, \varphi }(X_{\cO ,
\mu })$ and a commutative diagram
$$
\begin{matrix}
\pi _F: & X_{\cO , \mu } &           & \longrightarrow &       & \cH _F, \\
        &                & \stackrel{\pi _F}{\searrow  }&                 &  \stackrel{b_F}{\nearrow }     & \\
        &                &           & X_{\cO ,\mu }/W_0 &     &
\end{matrix}
$$
where $p_F: X_{\cO , \mu }\ra X_{\cO , \mu }/W_0$ is the natural
projection.

\end{enumerate}

We say that a closed point $[A, \iota , \cL ]$ in $X_{\cO , \mu }$ or
in any quotient of it is a {\em Heegner point} if $\End (A)\varsupsetneq \iota (\cO
)$. It was also shown in \cite{Ro2} that the morphisms
$b_F$ and $b_S$ are biregular on $X_{\cO , \mu }/W_0$ and $X_{\cO , \mu }/V_0(S)$,
respectively, outside a finite set of Heegner points.

It follows from these facts that the Galois
group $G=\Gal (k_{\cO }/k_{R_F})$ of the extension of
fields of moduli $k_{\cO }/k_{R_F}$ is naturally embedded in
$W_0$: any $\sigma \in G$ acts on a principally polarized abelian
variety with quaternionic multiplication
$(A, \iota :\cO \stackrel {\simeq}{\ra } \End
_{\qbar }(A), \cL )$ by leaving the $\qbar $-isomorphism class of
$\pi _F(A, \iota , \cL )=(A, \iota _{|R_F}: R_F\hookrightarrow
\End _{\qbar }(A), \cL )$ invariant.

Similarly, $\Gal (k_{\cO }/k_{S})$ embeds in $V_0(S)$ for any totally real order $S$
embedded in $\cO $. In what follows, we will describe the structure of the
groups $W_0$ and $V_0(S)$ attached to a polarized order $(\cO , \mu )$.
This will automatically yield Theorem \ref{main}. In fact, in Propositions
\ref{nontw} and \ref{tw}, we will be
able to conclude a rather more precise statement than the one given in
Section \ref{mm}.

The next proposition shows that the situation is simplified
considerably in the non-twisting case.

\begin{proposition}
\label{nontw}

Let $(A, \cL )$ be a principally polarized abelian variety over $\qbar
$ with quaternionic multiplication by $\cO $. Let $\iota : \cO \simeq \End (A)$
be any fixed isomorphism and let $\mu \in \cO $ be such that  $\mu
^2+D=0$ for some $D\in F^*_+$, $\disc (B) = (D)$, and $\iota (\beta )^*
= \iota (\mu ^{-1}\bar {\beta } \mu)$ for all $\beta \in \cO $.

If $(\cO , \mu)$ is a non twisting polarized order, then
$k_{\cO }=k_S$ for any totally real quadratic order $S\subset \cO $ over $R_F$
and

$$
\Gal (k_{\cO }/k_{R_F})\subseteq C_2^{\omega _{odd}}.
$$

\end{proposition}

{\em Proof. } It is clear from definition \ref{twist} that the
groups of twisting involutions $V_0(S)$  are trivial for any
subring $S$ of $\cO $. Since $\Gal (k_{\cO }/k_{S})\subseteq V_0(S)$, this
yields the first part of the proposition. As for the second,
since $\Gal (k_{\cO }/k_{R_F})\subseteq W_0$, we have that the
Galois group $\Gal (k_{\cO
}/k_{R_F})$ is contained in $U_0\subseteq W^1$, which is a
$2$-torsion abelian finite group. Our claim now follows from the
following lemma, which holds true for arbitrary pairs $(\cO , \mu
)$. $\Box $

\begin{lemma}
Let $(\cO , \mu)$ be a principally polarized maximal order. Then
$U_0\simeq C_2^{\omega _{odd}}$.
\end{lemma}

{\em Proof. } Let us identify $F(\mu )$ and $F(\sqrt {-D})$
through any fixed isomorphism. As $U_0$ naturally embeds in
$F(\sqrt {-D})^*/ F^* \Omega $, we first show that the maximal
$2$-torsion subgroup $H$ of $F(\sqrt {-D})^*/F^* \Omega $ is
isomorphic to $C_2^{\omega _{odd}}$.

If $\om \in F(\sqrt {-D})^*$ generates a subgroup of $F(\sqrt {-D})^*/F^* \Omega
$ of order $2$, then $\om ^2=\lambda \xi $ for some root of unity $\xi \in \Omega
$ and $\lambda \in F^*$. In particular, note that if $\om \in F(\sqrt
{-D})^*$, then $\om ^2\in F^*$ if and only if $\om \in F^*\cup F^*\sqrt
{-D}$. Let us write $\bar {H}=
H/\langle \sqrt {-D}\rangle $.

We then have that, if $\xi \in \Omega $, there exists at most a
single subgroup $\langle \om \rangle \subseteq \bar {H}$ such that
$\om \in F(\sqrt {-D})^*$, $\om ^2\in F^* \xi $. Indeed, if $\om
_1$, $\om _2\in F(\sqrt {-D})$, $\om _i^2=\lambda _i \xi $ for
some $\lambda _i\in F^*$ then $\frac {\om _1}{\om _2}^2\in F^*$
and hence $\frac {\om _1}{\om _2}\in F^*\cup F^* \sqrt {-D}$. This
shows that $[\om _1]=[\om _2]\in \bar {H}$.

Observe further that, if $\xi _f\in \Omega $ is a root of unity of
odd order $f\geq 3$, then $\om = \xi _f^{\frac {f+1}{2}}\in
F(\sqrt {-D})^*$ generates a $2$-torsion subgroup of $F(\sqrt
{-D})^*/F^* \Omega $ such that $\om ^2=\xi _f$.

It thus suffices to show that $\bar {\bar {H}}=H/\langle \sqrt
{-D}, \{ \xi _f ^{\frac {f+1}{2}}\} _{f\geq 3 \mbox{ odd}}\rangle
$ is trivial. Let $\om \in F(\sqrt {-D})^*$, $\om ^2=\lambda \xi
$, $\xi $ a root of unity of primitive order $f\geq 1$. If $f$ is
$2$ or odd we already know that the class $[\om ]\in \bar {\bar
{H}}$ is trivial. Further, it can exist no $\xi \in F(\sqrt {-D})$
of order $f=2^N$, $N\geq 2$, because otherwise $\xi ^{2^{N-1}}$
would be a square root of $-1$ and we would have that $F(\sqrt
{-D})=F(\sqrt {-1})$. This is a contradiction since $D R_F= \wp
_1\cdot ...\cdot \wp _{2 r}$, $r>0$.

Finally, it is also impossible that there should exist $\om \in
F(\sqrt {-D})$, $\om ^2=\lambda \xi $, $\xi ^f=1$, $f=2^N f_0$
with $N\geq 1$ and $f_0\geq 3$ odd. Indeed, since in this case
$\xi '=\xi ^{2^N}\in F(\sqrt {-D})$ is a primitive root of unity
of order $f_0$, $\om '=\xi ^{\frac {f_0+1}{2}}$ satisfies ${\om
'}^2=\xi '$. Then we would have $\frac {\om '}{\om }^2=(\xi
^{2^{N-1}}) \xi $ and this would mean that $[\frac {\om '}{\om
}]=[\om ]\in \bar {H}$, which is again a contradiction. This shows
that $\bar {\bar H}$ is trivial and therefore $H=\langle \sqrt
{-D}, \{ \xi _f ^{\frac {f+1}{2}}\} _{\xi _f\in \Omega
_{odd}}\rangle $. In order to conclude the lemma, we only need to
observe that both $\mu $ and $\xi _f ^{\frac {f+1}{2}}\in F(\mu )$
normalize the maximal order $\cO $ for any odd $f$, because their
respective reduced norms divide the discriminant $D$. $\Box $

\begin{corollary}
\label{quad}

Let $(\cO , \mu )$ be a non-twisting polarized order and
assume that $F(\sqrt {-D})$ is a CM-field with no purely imaginary roots of
unity. Then, for any real quadratic order $S$ over $R_F$,
$k_{\cO }/k_{R_F}=k_{S}/k_{R_F}$ is at most a quadratic extension.

If, in addition, $k_{R_F}$ admits a real embedding, then $k_{\cO }$
is a totally imaginary quadratic extension of $k_{R_F}$.

\end{corollary}

{\em Proof. } The first part follows directly from the above. As for
the second, it follows from a theorem of Shimura \cite{Sh1} which
asserts that the Shimura varieties $X_{\cO , \mu }$ fail to have real
points and hence the fields $k_{\cO }$ are purely imaginary. $\Box $

However, if on the other hand $(\cO , \mu )$ is twisting, the
situation is more subtle and less homogenous as we now show.

\begin{lemma}
\label{lemmatw}

Let $(\cO , \mu )$ be a {\em twisting} order in a totally
indefinite quaternion algebra $B$ over $F$ of discriminant $\disc
(B)=(D)$, $D\in F^*_+$. Then $U_0\subset V_0$ is a subgroup of
$V_0$ and $V_0/U_0\simeq U_0$. In particular, $W_0=V_0\simeq
C_2^{2 \om _{odd}}$.
\end{lemma}

{\em Proof. } Let $\om \in U_0$ be represented by an element $\om
\in \Norm _{F(\mu )^*}(\cO )\cap \cO $ and let $\nu \in V_0$ be a
twisting involution. We know that the class of $\nu $ in $\Norm
_{B^*}(\cO )/F^* \cO ^*$ is represented by a twist $\chi \in \Norm
_{B^*}(\cO )\cap \cO $ that satisfies $\chi ^2+\n (\chi )=0$ and
$\mu \chi =-\chi \mu $. Then we claim that $\om \nu \in V_0$ is
again a twisting involution of $(\cO , \mu )$. Indeed, first $\om
\chi \in \Norm _{B^*} (\cO )\cap \cO $, because both $\om $ and
$\chi $ do. Second, since $\om \in F(\mu )$, $\mu (\om \chi )= \mu
\om \chi = \om \mu \chi = -\om \chi \mu = -(\om \chi ) \mu $ and
finally, we have $\tr (\mu (\om \chi ))= \mu \om \chi + \bar {\om
\chi }\bar {\mu }= \mu \om \chi -\bar {\om \chi }\mu =-\tr {\om
\chi }\mu \in F$ and thus $\tr (\om \chi )=0$.

This produces a natural action of $U_0$ on the set of twisting
involutions of $(\cO , \mu )$ which is free simply because $B$ is
a division algebra. In order to show that it is transitive, let
$\chi _1$, $\chi _2$ be two twists. Then $\om =\chi _1 \chi
_2^{-1}\in F(\mu )$ because $\mu \om = \mu \chi _1 \chi
_2^{-1}=-\chi _1 \mu \chi _2^{-1}=\chi _1 \chi _2^{-1}\mu =\om \mu
$ and $F(\mu )$ is its own commutator subalgebra of $B$; further
$\om \in \Norm _{B^*}(\cO )$ because its reduced norm is supported
at the ramifying prime ideals $\wp \mid \disc (B)$. Let us remark
that, in the same way, $\chi _1 \chi _2\in \Norm _{F(\mu )^*}(\cO
)$.

We are now in a position to prove the lemma. Let $\nu \in V_0$ be
a fixed twisting involution. Then $U_0\subset V_0$: for any $\om
\in U_0$ we have already shown that $\om \nu $ is again a twisting
involution and hence $(\om \nu )\nu = \om \in V_0$ because $V_0$
is a $2$-torsion abelian group. In addition, the above discussion
shows that any element of $V_0$ either belongs to $U_0$ or is a
twisting involution and that there is a non-canonical isomorphism
$V_0/U_0\simeq U_0$. $\Box $

Observe that in the twisting case, by the above lemma, $U_0$ acts
freely and transitively on the set of twisting involutions of
$W^1$ with respect to $(\cO , \mu )$.

\begin{proposition}
\label{tw}

Let $(A, \cL )$ be a principally polarized abelian variety over $\qbar
$ with quaternionic multiplication by $\cO $. Let $\iota : \cO \simeq \End (A)$
be any fixed isomorphism and let $\mu \in \cO $ be such that  $\mu
^2+D=0$ for some $D\in F^*_+$, $\disc (B) = (D)$, and $\iota (\beta )^*
= \iota (\mu ^{-1}\bar {\beta } \mu)$ for all $\beta \in \cO $.

If $(\cO , \mu)$ is a twisting polarized order, let
$\chi _1$, ..., $\chi _{s_0}\in \cO $ be representatives of
the twists of $(\cO , \mu )$ up to multiplication by elements in $F^*$.
Then,

\begin{enumerate}
\item [(i)]
For any real quadratic order $S$, $S \not \subset F(\chi
_i)$, $1\leq i \leq s_0$,

$$
k_{\cO }=k_S
$$
\item[(ii)]
For any real quadratic order $S\subset F(\chi _i)\cap \cO $, $1\leq i \leq s_0$,
$k_{\cO }/k_{S_i}$ is (at most) a quadratic extension.

\item[(iii)]
$k_{\cO }= k_{S_1}\cdot ...\cdot k_{S_{s_0}}$ and
$\Gal (k_{\cO }/k_{R_F})\subseteq C_2^{2 \omega _{odd}}.$

\end{enumerate}

\end{proposition}

{\em Proof. } If $S\not \subset F(\chi _i)$ for any $i=1, ...,
s_0$, then $V_0(S)$ is trivial and hence,
since $\Gal (k_{\cO }/k_S)\subseteq V_0(S)$,
$\Gal (k_{\cO }/k_S)$ is also trivial. If, on the
other hand, $S\subseteq F(\chi _i)\cap \cO $, then
$V_0(S)\simeq C_2$ is generated by the
twisting involution associated to $\chi _i$. Again, we deduce that
in this case $k_{\cO }/k_S$ is at most a quadratic extension.

With regard to the last statement, note that $U_0\supseteq \langle
[\mu ]\rangle $ is at least of order $2$. Thus, if $(\cO , \mu )$
is a twisting polarized order, it follows from Lemma \ref{nontw}
that there exist two non-commuting twists $\chi $, $\chi '\in \cO
$. Then $R_F[\chi , \chi ']$ is a suborder of $\cO $ and, since
they both generate $B$ over $\Q $, the fields of moduli $k_{\cO }$
and $k_{R_F[\chi , \chi ']}$ are the same. This shows that $k_{\cO
}\subseteq k_{S_1}\cdot ...\cdot k_{S_{s_0}}$. The converse
inclusion is obvious.

Finally, we deduce that $k_{\cO }/k_{R_F}$ is a $(2, ...,
2)$-extension of degree at most $2^{2 \omega _{odd}}$ from Lemma
\ref{lemmatw}. $\Box $

\begin{remark} In the twisting case, the field of moduli of quaternionic
multiplication is already generated by the field of moduli of any
maximal real commutative multiplication but for finitely many
exceptional cases. This homogeny does not occur in the
non-twisting case.
\end{remark}

In view of corollaries \ref{nontw} and \ref{tw}, the shape of the
fields of moduli of the endomorphisms of the polarized abelian
variety $(A, \cL )$ differs considerably depending on whether it
gives rise to a twisting polarized order $(\cO , \mu )$ or not.

For a maximal order $\cO $ in a totally indefinite
quaternion algebra $B$ of principal reduced discriminant $D\in F^*_+$,
it is then obvious to ask the questions whether

\begin{enumerate}
\item[(i)] There exists $\mu \in \cO $, $\mu ^2+D=0$ such that
$(\cO , \mu )$ is twisting and,

\item[(ii)]
If $(\cO , \mu )$ is twisting, what is its twisting group $V_0$.
\end{enumerate}

Both questions are particular instances of the ones considered in
the appendix at the end of the article.

\section{Fields of moduli versus fields of definition}\label{defi}

In dimension $2$, the results of the previous sections are
particularly neat and can be made more complete. Let $C/\qbar $ be
a smooth irreducible curve of genus $2$ and let $(J(C), \Theta
_C)$ denote its principally polarized Jacobian variety. Assume
that $\End _{\qbar }(J(C))=\cO $ is a maximal order in an
(indefinite) quaternion algebra $B$ over $\Q $ of reduced
discriminant $D=p_1\cdot ... \cdot p_{2 r}$. Recall that $\cO $ is
unique up to conjugation or, equivalently by the Skolem-Noether
Theorem, up to isomorphism.

Attached to $(J(C), \Theta _C)$ is the polarized order $(\cO , \mu
)$, where $\mu =c_1(\Theta _C)\in \cO $ is a pure quaternion of
reduced norm $D$. As we have seen, a necessary condition for $(\cO
, \mu )$ to be twisting is that $B\simeq (\frac {-D,m}{\Q })$ for
some $m\mid D$. The isomorphism occurs if and only if for any odd ramified prime
$p \mid D$, $m$ is not a square mod $p$ if $p\nmid m$ (respectively $D/m$
if $p\mid m$).

In the rational case, the Atkin-Lehner and the positive
Atkin-Lehner groups coincide and $W=W^1=\{ \om _d; d\mid D\}
\simeq C_2^{2 r}$ is generated by elements $\om _d\in \cO $, $\n
(\om _d)=d\mid D$. Moreover, $U_0=\langle \om _D\rangle \simeq
C_2$.

If $(\cO , \mu )$ is a non twisting polarized order, then the
field of moduli of quaternionic multiplication $k_{\cO }$ is at
most a quadratic extension over the field of moduli $k_C$ of the
curve $C$ by Proposition \ref{nontw}. Moreover, $k_{\cO } = k_S$
for any real quadratic order $S\subset \cO $.

On the other hand, if $(\cO , \mu )$ is twisting and $B=(\frac
{-D, m}{\Q })$ for $m\mid D$, then $V_0=\{ 1, \om _m, \om
_{D/m},\om _D\} \simeq C_2^2$, where we can choose representatives
$\om _m$, $\om _{D/m}$ in $\cO $ such that $\mu \om _m=-\om _m \mu
$ and $\mu \om _{D/m}=-\om _{D/m} \mu $. Note that, up to
multiplication by non zero rational numbers, $\om _m$ and $\om
_{D/m}$ are the only twists of $(\cO , \mu )$. When we
particularize Proposition \ref{tw} to the case of Jacobian
varieties of curves of genus $2$, we obtain the following

\begin{theorem}

Let $C/\qbar $ be a smooth irreducible curve of genus $2$ such
that $\End _{\qbar }(J(C))=\cO $ is a maximal order in a rational
quaternion division algebra $B$ of reduced discriminant $D$.
Assume that the polarized order $(\cO, \mu )$ attached to $(J(C),
\Theta_C)$ is twisting and let $m\mid D$ be such that $B\simeq
(\frac{-D, m}{\Q })$. Then

\begin{enumerate}

\item [(i)]
$k_{\cO }/k_C$ is at most a quartic abelian extension.

\item [(ii)]
$k_{\cO }=k_S$ for any real quadratic order $S\not \subset \Q (\om
_m)\simeq \Q (\sqrt {m})$ or $\Q (\om _{D/m})\simeq \Q (\sqrt
{D/m})$.

\item [(iii)]
$k_{\Z [\om _m]}$ and $k_{\Z [\om _{D/m}]}$ are at most
quadratic extensions of $k_C$ and these are such that $k_{\cO
}=k_{\Z [\om _m]}\cdot k_{\Z [\om _{D/m}]}$.
\end{enumerate}

\end{theorem}

In \cite{Me}, Mestre studied the relation between the field of
moduli $k_C = k_{J(C), \Theta _C}$ of a curve of genus $2$ and its
possible fields of definition, under the sole hypothesis that the
hyperelliptic involution is the only automorphism on the curve.
Mestre constructed an obstruction in $\mathrm{Br }_2(k_C)$ for $C$
to be defined over its field of moduli. On identifying this
obstruction with a quaternion algebra $H_C$ over $k_C$, he showed
that $C$ admits a model over a number field $K$, $k_C\subseteq K$,
if and only if $H_C\otimes K\simeq \M _2(K)$.

If $\Aut (C)\varsupsetneq C_2$, Cardona \cite{CaQu} has recently
proved that $C$ always admits a model over its field of moduli
$k_C$.

Assume now, as in the theorem above, that $\End _{\qbar}(J(C))\simeq \cO $ is a
maximal order in a quaternion division algebra $B$
over $\Q $. Let $K$ be a field of definition of $C$; note that,
since $\End _{\qbar }(J(C))\otimes \Q = B$ is division, $\Aut
(C)\simeq C_2$ and therefore $k_C$ does not need to be a possible
field of definition of the curve. Having made the choice of a
model $C/K$, there is a minimal (Galois) field extension $L/K$ of
$K$ such that $\End _L(J(C))\simeq \cO $. This gives rise to a
diagram of Galois extensions

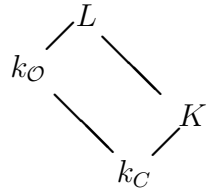
\vspace{0.5cm}
\begin{picture}(160,60)

\put(300,10){\makebox(0,0){\bfseries }}
\put(122,60){\makebox(0,0){$L$}}
\put(100,40){\makebox(0,0){$k_{\cO }$}}
\put(140, 0){\makebox(0,0){$k_C$}}
\put(162,22){\makebox(0,0){$K$}}

\thicklines
\put(128,52){\line(1,-1){22}}
\put(107,47){\line(1,1){10}}
\put(110,30){\line(1,-1){22}}
\put(147,7){\line(1,1){10}}
\end{picture}
\vspace{0.5cm}

The nature of the Galois extensions $L/K$ was studied in
\cite{DiRo} and \cite{DiRo2}, while the relation between the field
of moduli $k_{\cO }$ and the possible fields of definition $L$ of
the quaternionic multiplication was investigated by Jordan in
\cite{Jo}. In the next proposition we recall some of these facts,
and we prove that $L$ is the compositum of $K$ and the field of
moduli $k_{\cO }$.

\begin{proposition}\label{last}\footnote{Erratum: In the official published version, the statement of
Proposition 4.2 and its proof are incorrect. I heartily thank
Hakan Granath for pointing out to me the mistakes. In the present
version of the article, I have restated Proposition 4.2 and
provided a new proof for it.}

Let $C/K$ be a smooth curve of genus $2$ over a number field $K$
and assume that $\End _{\qbar }(J(C))$ is a maximal quaternionic
order $\cO $. Let $L/K$ the minimal extension of $K$ over which
all endomorphisms of $J(C)$ are defined. Then

\begin{enumerate}
\item[(i)] $\Gal (L/K)\simeq \{ 1\}$, $C_2$ or $D_2=C_2\times
C_2$.

\item[(ii)] $B\otimes _{\Q }L \simeq \mathrm{M}_2(L)$ and $L =
k_{\cO }\cdot K$.

\end{enumerate}

\end{proposition}

{\em Proof. } Statement $(i)$ was proved in \cite{DiRo}. As for
$(ii)$, let $M$ be any number field. Jordan proved in \cite{Jo}
that the pair $(J(C), \End (J(C)))$ admits a model over $M$ if and
only if $M$ contains $k_{\cO }$ and $M$ splits $B$. Since $A$ is
defined over $K$ and all its endomorphisms are defined over $L$,
we obtain that $L\supseteq k_{\cO }\cdot K$ and $B\otimes _{\Q }L
\simeq \mathrm{M}_2(L)$.

Let us now show that $L = k_{\cO }\cdot K$. By $(i)$, $\Gal
(L/k_{\cO }\cdot K)\subseteq \Gal (L/K)\subseteq D_2$. Assume on
the contrary that $L \varsupsetneq k_{\cO }\cdot K$; we will
encounter a contradiction. Let $\sigma \in \Gal (L/k_{\cO }\cdot
K)$, $\sigma \ne 1$, be such that $\sigma ^2=1$. Since $A$ is
defined over $K$, and according to the definition of $k_{\cO }$,
there exists an automorphism $\phi $ of the polarized abelian
variety $(A, \cL )$ such that $\alpha ^{\sigma }\cdot \phi = \phi
\cdot \alpha $ for any $\alpha \in \cO =\End_L(A)$. Let $L^{\sigma
}$ denote the fixed field of $L$ by $\sigma $. By \cite{DiRo},
Theorem 1.3, $\End _{L^{\sigma }}(A)\otimes \Q \simeq \Q
(\sqrt{m})$  for $m=-D$ or a positive divisor $m\mid D$, $m\ne 1$.
More precisely, as it is explained in Section 2 of \cite{DiRo},
$\End _{L^{\sigma }}(A)\otimes \Q = \Q (\gamma )$ for some $\gamma
\in B^*$, $\gamma ^2 = m$, such that $\alpha ^{\sigma } = \gamma
^{-1}\alpha \gamma $ for all $\alpha \in \cO $.

Hence, $\alpha ^{\sigma }\cdot \phi = \phi \cdot \alpha
\Leftrightarrow \gamma ^{-1}\alpha \gamma \cdot \phi = \phi \cdot
\alpha $. For $\alpha = \gamma $ we obtain that $\gamma \phi =
\phi \gamma $ and hence $\phi \in \Q (\gamma )\simeq \Q
(\sqrt{m})$.

By the Skolem-Noether theorem, there exists $\alpha \in \cO $ such
that $\alpha \gamma = -\gamma \alpha $. For such an element
$\alpha $, we obtain from the above equality that $\phi \gamma =
-\gamma \phi $. Since we also know that $\phi \in \Q (\gamma )$,
it follows that $\phi \in \Q \cdot \gamma $. This enters in
contradiction with the fact that $\phi $ is an automorphism. $\Box
$

\begin{example} Let $C$ be the smooth projective curve of hyperelliptic
model

$$
Y^2 = (1/48)X(9075X^4+3025(3+2 \sqrt{-3}) X^3 -6875 X^2 +220
(-3+ 2 \sqrt{-3}) X + 48).
$$

Let $A=J(C)/K$ be the Jacobian variety of $C$ over $K=\Q (\sqrt
{-3})$. By \cite{HaMu}, $A$ is an abelian surface with
quaternionic multiplication by a maximal order in the quaternion
algebra of discriminant $10$\footnote{In the published version,
there is a misprint: I wrongly claimed the discriminant to be
$6$.}. As it is explicitly shown in \cite{HaMu}, there is an
isomorphism between $C$ and the conjugated curve $C^{\tau}$ over
$\Q $. Hence, the field of moduli $k_C=\Q $ is the field of
rational numbers. By applying the algorithm proposed by Mestre in
\cite{Me}, we also obtain that the obstruction $H_C$ for $C$ to
admit a model over $\Q $ is not trivial. Hence $K=\Q (\sqrt{-3})$
is a minimal field of definition for $C$, though $C$ also admits a
model over any other quadratic field $K'$ that splits $H_C$.

In addition, it was shown in \cite{DiRo} that $L= \Q( \sqrt{-3},
\sqrt{-11})/K$ is the minimal field of definition of the
quaternionic endomorphisms of $A$. By a result of Shimura, Shimura
curves fail to have rational points over real fields. Hence,
$k_{\cO }$ must be a subfield of $L$ that does not admit a real
embedding. By Proposition \ref{last}, $L = k_{\cO }(\sqrt{-3})$
and thus $k_{\cO }$ is either $\Q (\sqrt{-11})$ or $L$ itself. It
would be interesting to determine which of these two fields is
$k_{\cO }$.\footnote{In the published version, I claimed that
$k_{\cO }=\Q (\sqrt{-11})$. This followed from the incorrect
statement of Proposition 4.2.}
\end{example}

\section{Appendix: Integral quaternion basis and distance ideals}

A quaternion algebra over a field $F$ is a central simple algebra
$B$ over $F$ of $\mathrm{rank} _F(B) = 4$. However, there are
several classical and more explicit ways to describe them which we
now review. Indeed, if $L$ is a quadratic separable algebra over
the field $F$ and $m\in F^*$ is any non zero element, then the
algebra $B = L + L e$ with $e^2=m$ and $e \beta = {\beta }^{\sigma
} e$ for any $\beta \in L$, where $\sigma $ denotes the
non-trivial involution on $L$, is a quaternion algebra over $F$.
The classical notation for it is $B = (L, m)$. Conversely, any
quaternion algebra over $F$ is of this form (\cite{Vi}).

In addition, if $\mathrm{char } (F)\not = 2$, then

$$
B = \left ( \frac {a, b}{F}\right ) = F + F i + F  j + F i j,
$$
with $i j = -j i$ and $i^2 = a\in F$, $j^2 = b\in F$ for any two
elements $a$, $b\in F^*$ is again a quaternion algebra over $F$
and again any quaternion algebra admits such a description. Note
that the constructions are related since $B=(\frac {a, b}{F})=
(F(i), b)$.

On a quaternion algebra $B$ there is a canonical anti-involution
$\beta \mapsto \bar {\beta }$ which is characterized by the fact
that, when restricted to any embedded quadratic subalgebra
$L\subset B$ over $F$, it coincides with the non-trivial
$F$-automorphism of $L$. Thus, if $B = (L, m)$, then $\bar {\beta
}=\overline {\beta _1 + \beta _2 e}= {\beta _1}^{\sigma }- {\beta
_2}^{\sigma } e$. The reduced trace and norm on $B$ are defined by
$\tr (\beta )=\beta +\bar {\beta }$ and $\n (\beta )=\beta \bar
{\beta }$.

Assume that $F$ is either a global or a local field of
char$(F)\not = 2$ and let it be the field of fractions of a
Dedekind domain $R_F$. An order $\cO $ in a quaternion algebra $B$
is an $R_F$-finitely generated subring such that $\cO \cdot F =
B$. Elements $\beta \in \cO $ are roots of the monic polynomial
$x^2-\tr (\beta ) x + \n (\beta )$, $\tr (\beta )$, $\n (\beta
)\in R_F$. We are now able to formulate the following question.

$\\ $ {\bfseries Question:} Let $B$ be a quaternion algebra over a
global or local field $F$, char $(F)\not =2$, and let $\cO $ be an
order in $B$.

\begin{enumerate}
\item
If $B\simeq (\frac {a, b}{F})$ for some $a, b\in R_F$, can one
find integral elements $\iota $, $\eta \in \cO $ such that $\iota
^2=a$, $\eta ^2=b$, $\iota \eta  = - \eta \iota $?
\item
If $B\simeq (L, m)$ for a quadratic separable algebra over $F$ and
$m\in R_F$, can one find $\chi \in \cO $ such that $\chi ^2=m$,
$\chi \beta = \bar {\beta } \chi $ for any $\beta \in L$?
\end{enumerate}
$\\ $
We note that Question 2 may be considered as a refinement of
Question 1. Indeed, let $\cO $ be an order in $B=(\frac {a,
b}{F})$ and fix an arbitrary element $i\in \cO $ such that
$i^2=a$. Then, while Question $1$ asks whether there exist
arbitrary elements $\iota $, $\eta \in \cO $ such that $\iota ^2=a$, $\eta ^2=b$ and
$\iota \eta = -\eta \iota $, Question $2$ wonders whether such an integral basis
exists with $\eta = i$.

If $B = \left ( \frac {a, b}{F}\right ) = F + F i + F  j + F i j$,
let $\cO _0 = R_F[i, j]$. Obviously, Question 1 is answered
positively whenever $\gamma ^{-1}\cO \gamma \supseteq \cO _0$ for
some $\gamma \in B^*$. The following proposition asserts that this
is actually a necessary condition. Although it is not stated in
this form in \cite{ChFr2}, it is due to Chinburg and Friedman, and
follows from the ideas therein. It is a consequence of Hilbert's
Satz 90. Let us agree to say that two orders $\cO $, $\cO '$ of
$B$ are of the same {\em type} if $\cO =\gamma ^{-1}\cO '\gamma $
for some $\gamma \in B^*$.

\begin{proposition}
\label{90}

Let $B = F + F i + F  j + F i j = \left ( \frac {a, b}{F}\right )
$ with $a$, $b\in R_F$. Let $\cO _0=R_F [i, j]$.

An order $\cO $ in $B$ contains a basis $\iota $, $\eta \in \cO $,
$\iota ^2=a$, $\eta ^2=b$, $\iota \eta  = - \eta \iota $ of $B$
if, and only if, the type of $\cO _0$ is contained in the type of
$\cO $.

\end{proposition}

{\em Proof. } Assume that there exist $\iota $, $\eta \in \cO $
satisfying the above relations. By the Skolem-Noether Theorem
(\cite{Vi}), $j$ and $\eta $ are conjugated (by, say, $\alpha \in
B^*$). Thus, by replacing $i$ by $\alpha ^{-1}i \alpha $ and $\cO
_0$ by $\alpha ^{-1}\cO _0\alpha $, we may assume that $j=\eta \in
\cO $. We then need to show the existence of an element $\gamma
\in F(j)=F(\eta )$ such that $\gamma ^{-1}\iota \gamma = i$.

We have $i \eta = -\eta i$ and thus $\eta = -i^{-1} \eta i$. In
addition, since $\iota \eta = -\eta \iota $, $\iota i^{-1} \eta i
= \eta \iota $. Hence, $(\iota i^{-1}) \eta = \eta (\iota i^{-1})$
and we deduce that $\iota i^{-1}\in F(\eta )$ is an element of
norm $\Norm _{F(\eta )/F}(\iota i^{-1})=1$.

By Hilbert's Satz 90, there exists $\om \in F(\eta )$ such that
$\iota i^{-1} = \om \bar {\om }^{-1}$, that is, $\iota = \om \bar
{\om }^{-1} i$. Stated in this form, we need to find an element
$\gamma \in F(\eta )$ with $\gamma ^{-1}\om \bar {\om }^{-1} i
\gamma = i$. Since $\gamma i = i\bar {\gamma }$, we can choose
$\gamma =\om $. $\Box $

An order $\cO $ in $B$ is {\em maximal} if it is not properly
contained in any other. It is an {\em Eichler order} if it is the
intersection of two maximal orders. The reduced discriminant ideal
of an Eichler order is $\disc (\cO )=\disc (B)\cdot \cN $ for some
integral ideal $\cN $ of $F$, the {\em level} of $\cO $, coprime
to $\disc (B)$ (see \cite{Vi}, p. 39). With this notation, maximal
orders are Eichler orders of level $1$.

\begin{corollary}

Assume that $F$ is a local field and that $\cO $ is an Eichler
order of level $\cN $ in $B=(\frac {a, b}{F})$, $a$, $b\in R_F$.
Then, there exist $\iota $, $\eta \in \cO $, $\iota ^2=a$, $\eta
^2=b$, $\iota \eta = -\eta \iota $ if and only if $\cN \mid 4 a
b$.

\end{corollary}

{\em Proof. } By \cite{Vi}, $\S 2$, there is only one type of
Eichler orders of fixed level $\cN $ in $B$. Remark that, if $B$
is division, necessarily $\cN =1$. Let $\cO _0=R_F [i, j]$. Since
$\disc (\cO _0)= 4 a b$, as one can check, a necessary and
sufficient condition on $\cO $ to contain a conjugate order of
$\cO _0$ is that $\cN \mid 4 a b$. The corollary follows from
proposition \ref{90}. $\Box $

In the global case, the approach to Question 1 can be made more
effective under the assumption that $B$ satisfies the Eichler
condition. Namely, suppose that some archimedean place $v$ of $F$
does not ramify in $B$, that is, $B\otimes_F F_v\simeq \M
_2(F_v)$. Here, we let $F_v\simeq \R $ or $\C $ denote the
completion of $F$ at $v$.

The following theorem of Eichler describes the set $\cT (\cN )$ of
types of Eichler orders of given level $\cN $ purely in terms of
the arithmetic of $F$. Let $\Pic _+(F)$ be the narrow class group
of $F$ of fractional ideals up to principal fractional ideals
$(a)$ generated by elements $a\in F^*$ such that $a>0 $ at any
real archimedean place $v$ that ramifies in $B$ and let
$h_+(F)=|\Pic _+(F)|$.

\begin{definition}

The group $\overline {\Pic }_+^{\cN} (F)$ is the quotient of $\Pic
_+(F)$ by the subgroup generated by the squares of fractional
ideals of $F$, the prime ideals $\wp $ that ramify in $B$ and the
prime ideals $\mathfrak q$ such that $\cN $ has odd $\mathfrak
q$-valuation.

\end{definition}

The group $\overline {\Pic }_+^{\cN} (F)$ is a $2$-torsion finite
abelian group. Therefore, if $h_+(F)$ is odd, then $\overline
{\Pic }_+^{\cN} (F)$ is trivial.

\begin{proposition}[\cite{Ei1}, \cite{Ei2}, \cite{Vi}, p. 89]
\label{type} The reduced norm $\n $ induces a bijection of
sets

$$
\cT (\cN )\stackrel {\sim }{\longrightarrow }\overline {\Pic
}_+^{\cN }(F).
$$

\end{proposition}

The bijection is not canonical in the sense that it depends on the
choice of an arbitrary Eichler order $\cO $ in $B$. For $\cN =1$,
the bijection is explicitly described as follows. For any two
maximal orders $\cO $, $\cO '$ of $B$ over $R_F$, define the {\em
distance ideal} $\rho (\cO , \cO ')$ to be the order-ideal of the
finite $R_F$-module $\cO /\cO \cap \cO '$ (\cite{Re}, p. 49).
Alternatively, $\rho (\cO , \cO ')$ can also be defined locally in
terms of the local distances between $\cO \otimes _{R_F}R_{F_{\wp
}}$ and $\cO '\otimes _{R_F}R_{F_{\wp }}$ in the Bruhat-Tits tree
$\cT _{\wp }$ for any (non-archimedean) prime ideal $\wp $ of $F$
that does not ramify in $B$ (\cite{ChFr1}). Finally, $\rho (\cO ,
\cO ')$ is also the {\em level} of the Eichler order $\cO \cap \cO
'$. This notion of distance proves to be suitable to classify the
set of types of maximal orders of $B$, as the assignation $\cO
'\mapsto \rho (\cO , \cO ')$ induces the bijection claimed in
proposition \ref{type}.

\begin{corollary}

Let $B = \left ( \frac {a, b}{F}\right )$, $a$, $b\in R_F$ be a
quaternion algebra over a global field $F$. If $B$ satisfies
Eichler's condition and $h_{+}(F)$ is odd then, for any Eichler
order $\cO $ in $B$, there is an integral basis $\iota $, $\eta
\in \cO $, $\iota ^2=a$, $\eta ^2=b$, $\iota \eta = -\eta \iota $
of $B$.

\end{corollary}

As for Question $2$, let $B = F+F i +F j +F i j= \left ( \frac {a,
b}{F}\right )= (L, b)$ with $a$, $b\in R_F$ and $L=F(\sqrt {a})$.
Choose an arbitrary order $\cO $ of $B$. For given $\eta \in \cO
$, $\eta ^2=a$, we ask whether there exists $\chi \in \cO $, $\chi
^2=b$, such that $\eta \chi = -\eta \chi $. By proposition
\ref{90}, a necessary condition is that $\cO _0 = R_F[i,
j]\subseteq \cO $ up to conjugation by elements of $B^*$ and,
without loss of generality, we assume that this is the case. With
these notations, we have

\begin{definition}

Let $\cO \supseteq \cO '$ be two arbitrary orders in $B$. The
transportator of $\cO '$ into $\cO $ over $B^*$ is $(\cO :\cO ')
:=\{ \gamma \in B^*, \gamma ^{-1}\cO _0 \gamma \subset \cO \} $.

\end{definition}

Note that $\Norm _{B^*}(\cO )$ is a subgroup of finite index
of $(\cO :\cO ')$.

\begin{proposition}
\label{Q2}

Let $\cO \supseteq \cO _0$ be an order in $B$ and let $\eta \in
\cO $, $\eta ^2=a$. Then, there exists $\chi \in \cO $, $\chi ^2=b$, $\eta \chi =
-\chi \eta $ if and only if $\eta = \gamma ^{-1} i \gamma $ for $\gamma \in (\cO
:_{B^*}\cO _0)$.

\end{proposition}

Let $f=|(\cO :\cO _0) : \Norm _{B^*}(\cO )|$
be the index of the normalizer group $\Norm _{B^*}(\cO )$ in $(\cO
:\cO _0)$. Let $\mathcal E (a)$ be the finite set of
$\Norm _{B^*}(\cO )$-conjugation classes of elements $\eta \in \cO $ such
that $\eta ^2=a$. Then, it follows from the above proposition that Question
$2$ for $(\cO , \eta )$
is answered in the affirmative for elements $\eta $ lying on exactly $f$
of the conjugation classes in $\mathcal E (a)$.
Again, the cardinality of $\mathcal E (a)$ can be explicitly
computed in many cases in terms of class numbers by means of
the theory of Eichler optimal embeddings (cf. \cite{Vi}).

{\bfseries Acknowledgements.} I am indebted to P. Bayer for her
assistance throughout the elaboration of this work. I also express
my gratitude to E. Friedman, J. Brzezinski, H. Granath and A.
Arenas for some helpful conversations. Finally, I thank J. Kramer and
U. Kuehn for their
warm hospitality at the Humboldt-Universit\"{a}t zu Berlin during the
fall of 2001.

\end{document}